\documentclass[12pt,a4paper]{amsart}
\usepackage{amssymb, amstext, amscd, amsmath, color}
\usepackage{url}
\usepackage{hhline}

\usepackage{amsfonts}
\usepackage{amsthm}
\usepackage{amsfonts}
\usepackage{array}
\usepackage{epsfig}
\usepackage{eucal}
\usepackage{latexsym}
\usepackage{mathrsfs}
\usepackage{textcomp}
\usepackage{verbatim}
\newcommand{\bR}{{\mathbb{R}}}

\newcommand{\bZ}{{\mathbb{Z}}}

  \newcommand{\G}{{\mathcal{G}}}
\renewcommand{\H}{{\mathcal{H}}}

\renewcommand{\L}{{\mathcal{L}}}

\renewcommand{\P}{{\mathcal{P}}}

\renewcommand{\S}{{\mathcal{S}}}
  \newcommand{\T}{{\mathcal{T}}}

\newcommand{\ul}{\underline  }
\newcommand{\coker}{\operatorname{coker}}
\newcommand{\trace}{\operatorname{trace}}

\newtheorem{thm}{Theorem}[section]

\newtheorem{cor}[thm]{Corollary}

\newtheorem{definition}[thm]{Definition}

\newtheorem{lem}[thm]{Lemma}

\newtheorem{prop}[thm]{Proposition}
\newtheorem{rem}[thm]{Remark}

\numberwithin{equation}{section}

\begin{document}

\markboth{J. C. Owen and S. C. Power}
{Frameworks, symmetry and rigidity}

%   \catchline

\title{FRAMEWORKS SYMMETRY AND RIGIDITY}
\footnote{Supported by a London Mathematical Society Scheme 7
Grant.}

\author{J. C. Owen}

\address{D-Cubed, Siemens PLM Software, Park House, \\
Castle Park, Cambridge, United Kindom\\
owen.john.ext@siemens.com}

\author{S. C. Power}

\address{Department of Mathematics and Statistics,\\Lancaster University,\\
Lancaster, LA1 4YF, United Kingdom \\
s.power@lancaster.ac.uk
}

\maketitle

%\pub{Received (received date)}{Revised (revised date)}
%{Communicated by (Name)}

\begin{abstract}
Symmetry equations are obtained for the rigidity matrix of a
bar-joint framework in $\bR^d$. These form the basis for a short proof  of
the Fowler-Guest symmetry group generalisation of the
Calladine-Maxwell counting rules.  Similar symmetry equations are
obtained for the Jacobian of diverse framework systems, including
constrained point-line systems that appear in CAD, body-pin
frameworks, hybrid systems of distance constrained objects
and infinite bar-joint frameworks. This leads to
generalised forms of the Fowler-Guest character formula
together with counting rules in terms of counts of
symmetry-fixed elements. Necessary conditions for
isostaticity are obtained for asymmetric frameworks, both
when symmetries are present in subframeworks
and when symmetries occur in partition-derived frameworks.

\end{abstract}

\section{Introduction}	%) A SECTION HEADING

Let $(G,p)$ be a \emph{framework} in $\mathbb{R}^{d}$ which, by
definition, consists of an abstract graph $G = (V,E)$
%, with $v$ vertices
%and $e$ edges, and distinct
and a vector $p=(p_{1}, \dots, p_{v})$ composed of framework
points in $\mathbb{R}^{d}$.
%framework points $p_{1}, \dots, p_{v}$ in $\mathbb{R}^{d}$. The
%framework points form the framework vector $p= (p_{1}, \dots,
%p_{v})$ in $\prod_{k=1}^{v} \mathbb{R}^{d}$.
When $(G,p)$ is viewed in the natural way as a pin-jointed bar
framework  in $\mathbb{R}^{d}$
%, and if $p$ is  generic,
then there is a counting condition for bars and joints that the
framework must satisfy if it is known to be isostatic, which is to
say that the structure is rigid in a natural sense
(infinitesimally rigid) and at the same time is not
overconstrained. More generally, in the nonisostatic case, there
is a single condition relating the four quantities, $v=|V|,
e=|E|$, the number $m$ of non-trivial independent infinitesimal
motions (also known as mechanisms), and the number $s$ of
independent stresses that the structure can carry.
For
$d=2$ this
is the extended Maxwell rule (Calladine\cite{cal})
\begin{equation}\label{e:cal}
m-s= 2v-e-3
\end{equation}
while for $d=3$ one has
%\begin{equation}\label{e:cal}
$m-s= 3v-e-6$.
%\end{equation}
The equations arise from a consideration of the kernel and
cokernel of the rigidity matrix for the framework and their
respective dimensions, $m$ and $s$.

 Recently, in the context of the analysis of loads and stresses
 in symmetric structures, Fowler and Guest
\cite{fow-gue} have obtained an extended counting rule for
symmetric frameworks in two and three dimensions and these formulae are
a source of additional necessary counting conditions.
%which are generic modulo a specific symmetry
%group ${\G}$. This
In three dimensions the formula takes the elegant form
\begin{equation}\label{e:FGequation}
\Gamma(m) - \Gamma(s)= \Gamma(v) \times \Gamma _{xyz}- \Gamma(e)-
\Gamma_{xyz}- \Gamma_{R_{x}R_{y}R_{z}}
\end{equation}
where each $\Gamma$ denotes a character list for a representation of
the rigid motion symmetry group ${\G}$ of the framework. Thus the
equation represents a set of equations, one for each element of
$\G$. The list $\Gamma(e)$, for example, arises from an elementary permutation
representation $\rho_e$ of $\G$ on a real vector space with basis indexed by the edges
of $G$. Specifically
$$\Gamma(e) = \trace(\rho_e(g_1), \dots , \trace(\rho_e(g_r))$$
for some choice of elements $g_1, \dots ,g_r$ of $\G$, typically a set of generating  elements with $g_1$ the identity element.

The significance of the formulae lie in the fact that the
right-hand side is readily computable depending only on the
abstract graph $G$ of the framework rather than the metrical
detail. In particular  $\trace(\rho_e(g_k))$
is the number of edges that are left unmoved by the
symmetry $g_k$.
The left hand side however carries information on the
possibilities for stresses and flexes. Evaluating the formula for
the identity element $g_1$ of $\G$  gives the Calladine-Maxwell rules.
See also Ceulemans and Fowler \cite{ceu-fow} for an analogous symmetry variant of
Euler's formula for polyhedra.

Our first purpose is to obtain an explicit symmetry
equation
\begin{equation}\label{e:symmequation}
R = \rho_e(g^{-1})R\hat{\rho}_v(g),\quad g \in \G,
\end{equation}
for the rigidity matrix $R=R(G,p)$ of a bar-joint framework in $\bR^d$,
which shows how the matrix intertwines
representations of $\G$ associated with the edges and with the vertices.
Here  $\hat{\rho}_v$ is the representation $\rho_n \otimes \rho_{sp}$ where $\rho_n$
is the natural permutation representation associated with the vertices (nodes)
and $\rho_{sp}$ is the usual orthogonal representation of $\G$ in $\bR^d$.
%For $d=3$ it provides
%the character list $\Gamma(v) \times \Gamma_{xyz}$.
From this we obtain a simple proof of a general Fowler-Guest
formula for frameworks in $\bR^d$. The proof is
coordinate free and in fact the unitary equivalence of subrepresentations
that underlies the formula may be implemented by the partially
isometric part of the polar decomposition $R = U(R^*R)^{1/2}$.

Our second purpose is to show that the method is versatile and
readily applicable to higher order frameworks. For example, we
consider body-bar frameworks
%such as the Stewart platform
%\cite{das-mru}
and constraint systems for geometric objects, such
as the constraints of geometries arising in CAD. Once again we
obtain symmetry equations, equivalent representations, character
formulae and counting conditions.

\begin{center}
\begin{figure}[h]
\centering
\includegraphics[width=8cm]{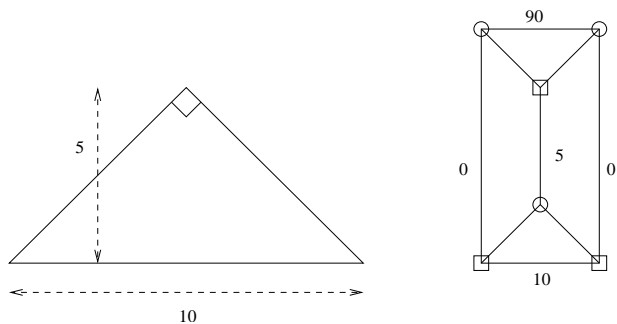}
\caption{The geometric drawing has an abstract graph in which
lines are represented by circular vertices and the points by
square vertices. The labeled edges represent angular dimensions of
90 degrees and distances of 0, 5 and 10. A count for the
reflection symmetry of the graph implies the singularity of the
equation system for the drawing.}
\end{figure}
\end{center}

Figure 1 shows a practical application in CAD for the symmetry adapted
Maxwell counting rule. The geometric figure on the right shows a
triangle which has reflection symmetry about a vertical axis. The graph
on the right shows the corresponding constraint graph taken from a CAD
constraint solving application. In this graph the square nodes represent
points, the circular nodes represent lines and the edges represent
either a distance or angle dimension with the specified value or a
point-line coincidence. Notice that this graph has a corresponding
two-fold symmetry. The equations represented by the geometric figure and
the constraint graph are singular for the following reason. If the
dimension with value 5 which specifies the height of the triangle is
removed then the apex of the triangle can be placed anywhere on a circle
with the base of the triangle as diameter (due to the perpendicular
constraint shown) and so this circle has radius 5. Thus in the symmetric
configuration shown the height dimension has attained its maximum
possible value and cannot be increased. This indicates singularity in
the equations. In this example $v=6$ and $e=9$ so the right hand side of
Equation (1.1) evaluates to zero and is compatible with $m=s=0$. For this
reflection symmetry we will show that the symmetry adapted equation takes the form
$g(m)-g(s)=g(e)-1$. Since $g(e)$, the number of edges of the graph which are
unchanged by the two-fold symmetry, is three, this equation requires
that $m>0$ which says that the equations have at least one infinitesimal
motion which means here that they are singular.

We also show how one may
obtain symmetry equations and the character formula
for infinite frameworks. In particular, in the case of periodic
frameworks we obtain
the periodic form
\begin{equation}\label{e:perFG}
 \Gamma_{p}(m) -\Gamma_{p}(s)=\Gamma_{p}(v)\cdot \Gamma(sp) - \Gamma_{p}(e)
 -\Gamma_{p}(rig),
\end{equation}
in which the trace lists are  associated with  finite-dimensional
representations of a translation subgroup quotient of the spatial symmetry group.
%Here too appropriately isostatic frameworks
%lead to the equivalence of the representations of the symmetry
%group that are associated with framework edges and vertices.

Finally we indicate how the symmetry analysis may be exploited
further in two distinct ways, and even for \textit{asymmetric}
frameworks. In the first we consider symmetries in vertex induced
subframeworks while in the second we consider latent symmetries in
partition-derived frameworks. For the symmetry group identity
element the properties of sub-frameworks and derived frameworks
both give the same  well-known necessary requirement for non-singularity of
the relevant Jacobian (such as, in two dimensions,  $2v-e\geq 3$ for every sub-graph with
$e$ edges and $v$ vertices). However, symmetry in subframeworks or
partition derived frameworks both give new and useful predictions.

For example we obtain in Theorem \ref{t:partitiontheorem}
a  "singularity predictor",
in the form of a set of necessary counting conditions for an isostatic
framework $(G,p)$ in $\bR^d$.
%CHECK LATER ALSO - SAY pi are distinct.
%\textit{If $(G,p)$ is proper (in that $p$ spans $\bR^d$) then for
% each proper subframework $(X,p)$ and each spatial symmetry $g$ of
%$(X,p)$},
%\[
%|trace(g).v_X^g - e_X^g - trace(\rho_{rig}(g)))| \leq dv_X-e_X
%-d(d+1)/2
%\]
%\textit{where $v_X^g$ (resp. $e_X^g$) is the number of vertices
%(resp. edges) in the graph $X$ that are unmoved by the symmetry}
%The quantity $trace(\rho_{rig}(g)))$ is framework independent and
%is derived from the representation of $g$ in the
%$d(d+1)/2$-dimensional space of rigid motions.
As a simple corollary of this we observe that
for a planar isostatic framework with a reflection symmetry $g$ in a subframework  $X$
%, for example the (framework
%independent) quantity  $trace(\rho_{rig}(g))$ is $-1$ and so in
%this case we obtain the useful isostaticity condition
we have the necessary condition
\[
|-e_X^g+1| \leq 2v_X - e_X-3.
\]
where $e_X^g$  is the number of edges of $X$ left unmoved by $g$.
%for all subframeworks $(X,p)$ and their reflection symmetries $g$.

Many authors have considered group representations in the analysis
of symmetric structures, often adopting symmetry adapted
coordinate spaces for stresses and flexes. See, for example,
Kangwai and Guest\cite{kan-gue} and the survey
Kangwai, Guest and Pellegrino\cite{kan-gue-pel}.  In this vein irreducible group  considerations were
introduced in the detailed engineering calculations of Kangwai and Guest\cite{kan-gue} and subsequently put into the useful character equational form by Fowler and Guest\cite{fow-gue}.
In contrast to this we bring out the symmetry equations as the essential
feature of symmetric bar-joint frameworks and we use them to identify invariant subspaces and thereby obtain a short derivation of the formula.
We note that Schulze\cite{sch-2} has given another rigourous proof
of the Fowler-Guest formula which uses a more expansive analysis of subspaces for the block diagonalisation of the rigidity matrix. Moreover interesting applications
are given to noninjective frameworks with coincident vertices which are not considered
here.

The usefulness of the Fowler-Guest formula has been shown in
Connelly, Fowler, Guest, Schulze and Whiteley\cite{con-et-al} where it was used
to derive a complete list of the necessary counting conditions
for bar-joint frameworks in two and three dimensions. These conditions are in terms of counts for the number of vertices or edges that are left unmoved by various symmetries.
In our Corollary \ref{countingCor} we recover some of these results while Theorems \ref{t:pointline}, \ref{t:body}, \ref{t:partitiontheorem} and \ref{t:last} lead to analogous counting constraints.

For a planar isostatic framework one has $m=s=0$ and
hence the necessary equality $2v-e-3=0$.  This is not a sufficient
condition as one also needs subframeworks not to be
overconstrained. However, it is a fundamental and celebrated
theorem of Laman\cite{lam} that the necessary count condition
$2v-e=3$ together with the  inequality $2v_X - e_X\geq 3$ for all
subgraphs $X$ is a sufficient condition for a
\textit{generic} framework to be isostatic.  Thus necessary and
sufficient conditions are known for the two dimensional generic
case. We do not consider sufficiency conditions below but we note that
 Schulze\cite{sch-3} has recently obtained
Laman theorems for frameworks in the plane with various symmetry.

For further background on rigidity and diverse constraint problems
see, for example, Asimow and Roth\cite{asi-rot}, Connelly et al\cite{con-et-al},
Graver, Servatius and Servatius\cite{gra-ser-ser}, Jackson and Jordan\cite{jac-jor}, Owen\cite{Owen}, Owen and Power\cite{owe-pow-1},
%Roth\cite{roth},
and Whiteley\cite{whi-2}.

 We would like to thank Simon Guest and Nadia Mazza for interesting
discussions and the anonymous referees for helpful comments.

\section{Frameworks and Symmetries.}

We begin with a formal introduction to mathematical bar-joint frameworks
$(G,p)$ in $\bR^d$, to the rigidity matrix $R(G,p)$ and to the
spatial symmetry group $\G$ of a framework. Also, viewing $\G$ as an abstract group we
consider elementary representations of $\G$ as permutation transformations
of vector spaces associated with the vertices and
with the edges.

\subsection{The rigidity matrix}
Let $G=(V,E)$, $n=|V|, m=|E|$ be a finite connected graph, with no
multiple edges. A \textit{framework} in $\mathbb{R}^{2}$ is a pair
$(G,p)$ where $p=(p_{1}, \dots, p_{n})$ is a framework vector with
framework points $p_{i}=(x_i,y_i)$ in $\mathbb{R}^{2}$ that are
associated with an ordering $v_1, \dots ,v_n$ of the vertices.
Thus we allow framework points to coincide. The \textit{rigidity
matrix} $R = R(G,p)$ for the framework $(G,p)$ is an $m \times 2n$
real matrix whose columns are labeled by $x_1, y_1, x_2, y_2,
\dots , x_n, y_n,$ and whose rows are labeled by some ordering
$e_1, \dots ,e_m$ of the edges. If $e= (v_{i}, v_{j})$ is an edge
of $G$ then the matrix entries of $R$ in the row for $e$ are zero
except possibly in the columns for $x_{i}, y_{i},x_{j}, y_{j}$
where we have, respectively, $x_{i}-x_{j}, y_{i}-y_{j},
x_{j}-x_{i}, y_{j}-y_{i}, 1 \leq i \leq n$. Thus for notational economy
we allow framework point coordinates to agree with
their labels.

The rigidity matrix gives a linear transformation from the $2n$-dimensional
real vector space
\[
{\H}_{v}= \sum_{k=1}^{n}\oplus (\mathbb{R}_{x_{k}} \oplus
\mathbb{R}_{y_{k}}),
\]
associated with the vertices,
 to the $m$-dimensional real vector space,
\[
 {\H}_{e}= \sum_{k=1}^{m} \oplus
\mathbb{R}_{e_{k}}
\]
associated with edges.  Here each vector space summand
$\mathbb{R}_{x_{k}}, \mathbb{R}_{y_{k}}, \mathbb{R}_{e_{k}}$ is a
copy of $\mathbb{R}$.
Let
 $\xi_{x_{k}}, \xi_{y_{k}},
1 \leq k \leq n$, denote the standard basis for $\H_v$ and write  $\xi_{e_{k}}, 1
\leq k \leq m,$ for the standard basis  for $\H_e$.
Then the matrix entry $x_i-x_j$ in row
$e=(v_i,v_j)$  and column $x_i$ is given by the standard inner product $\langle
R\xi_{x_i}, \xi_e \rangle$.

The rigidity matrix $R(G,p)$ of a framework $(G,p)$ in $\bR^d$ is defined
in exactly the same manner. Alternatively it may be defined
as one half of the Jacobian derivative of the nonlinear map
from $\H_v$ to $\H_e$ which is determined by the quadratic
distance equations for the framework. We adopt this viewpoint in
Section 4.

The rigidity matrix derives its name from the fact that vectors $u=(u_i)=(u_{x_i},u_{y_i})$ in its kernel (nullspace) are infinitesimal flexes in the following sense. They indicate
directions (or velocity directions) in which
for each edge the disturbances of edge length
\[
|p_i-p_j|-|(p_i+tu_{x_i}-(p_j+tu_{y_i})|
\]
is $O(t^2)$ as $t$ tends to zero.
Also, vectors in the cokernel (the kernel of the transpose matrix) correspond to self stresses.
Moreover we have the following fundamental definition.

\begin{definition}
A  framework $(G,p)$ in the plane (resp. in $\bR^3$) with graph
 $G=(V,E)$ is infinitesimally rigid if
the rank of $R(G,p)$ is  $2|V|-3$ (resp. $3|V|-6$) and is isostatic if it is
infinitesimally rigid and the rank of $R(G,p)$ is $|E|$.
\end{definition}

As an illustration we shall  keep in view the symmetric framework $(G,p)$ in $\bR^2$
indicated in Figure 1, with framework vector
\[
p=((2,0),(3,1),(4,0),(3,-1),(-4,0),(-3,1),(-2,0),(-3,-1)).
\]
The subframework on the first four vertices is
infinitesimally rigid as is its mirror image in the $y$-axis.
The entire framework appears to have one non-trivial
infinitesimal flex in addition to the three spatial flexes
and this is readily confirmed.
\begin{center}
\begin{figure}[h]
\centering
\includegraphics[width=9cm]{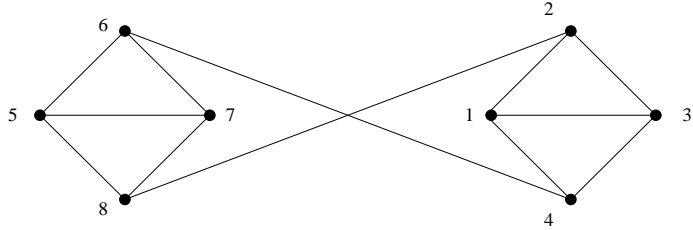}
\caption{A symmetric bar-joint framework.}
\end{figure}
\end{center}
The rigidity matrix has the form
$$R=\begin{bmatrix}
R_1&0\\0& R_1\\
T_1&T_2
\end{bmatrix}
$$
where $R_1$ is the $5$ by $8$ matrix
$$R_1=\begin{bmatrix}
-1&-1&1&1&0&0&0&0\\
0&0&-1&1&0&0&1&-1\\
0&0&0&0&1&1&-1&-1\\
-1&1&0&0&0&0&1&-1\\
2&0&0&0&-2&0&0&0
\end{bmatrix}
$$
and where the submatrix $\begin{bmatrix}
T_1& T_2
\end{bmatrix}$ is the $2$ by $16$ matrix corresponding to
the two long framework
edges $[p_2,p_8], [p_4,p_6]$.
%, has the
%form DOESN"T WORK ?
%$$\begin{bmatrix}
%0&0&6&2&0&0&0&0&0&0&0&0&0&0&-6&-2\\
%0&0&0&0&0&0&6&-2&0&0&-6&2&0&0&0&0
%\end{bmatrix}
%$$

\subsection{Graph symmetry.}
Let $G$ have vertices $v_1, \dots ,v_n$
and let $\sigma$ be
a permutation of $(1, \dots , n)$
corresponding to
an automorphism of $G$. We also write $\sigma : V\to V$ and
$\sigma :E\to E$ for the
corresponding bijective maps so that $\sigma(v_i) = v_{\sigma(i)}$.
Let $\sigma_e$
denote the associated linear transformation of $\H_e$,  where
$\sigma_e \xi_f = \xi_{\sigma (f)}$, and let $\sigma_e$ also denote
its representing matrix. The transformation and matrix
 $\sigma_v$ is similarly defined on the space
 $\H_v$  by the specification  $\sigma_v \xi_{x_i} = \xi_{x_{\sigma (i)}},
 \sigma_v \xi_{y_i} = \xi_{y_{\sigma (i)}}, 1 \le i \le n$.

We first note how $R(G,p)$ is transformed, even in the absence of
framework symmetry, on replacing the framework vector $p=(p_{1}, \dots, p_{n})$ by
$\sigma(p)=(p_{\sigma(1)}, \dots, p_{\sigma(n)})$.

\begin{lem}\label{l:symm} Let $(G,p)$ be a framework in $\bR^d$ with
 rigidity matrix $R(G, p)$ and let $\sigma$ be a graph
 automorphism. Then
\begin{equation}\label{e:graphsymmetry}
R(G,\sigma(p)) = \sigma_e^{-1} R(G,p) \sigma_v
\end{equation}
\end{lem}

\begin{proof}
For notational simplicity let $d=2$. Let $\sigma(p) =
(p_{\sigma(1)},\dots , p_{\sigma(n)}) = (p_1', \dots ,p_n'),$ and $
p_i' = (x_i',y_i'), 1 \leq i\leq n$. Associated with $e =
(v_i,v_j)$ we have $x_i'-x_j' = x_{\sigma(i)}- x_{\sigma(j)}$.
This difference appears in the $\sigma(e)$ row and $\sigma(x_i)$
column of $R(G,p)$ and so
\[
x_i'-x_j' =  \langle R(G,p)\xi_{\sigma(x_i)}, \xi_{\sigma(e)}
\rangle .
\]
On the other hand, from the definition of $R(G, \sigma(p))$,
\begin{align*}
x_i'-x_j' & =\langle R(G,\sigma(p))\xi_{x_i}, \xi_{e}
\rangle \\
&=  \langle R(G,\sigma(p))\sigma_v^{-1}\xi_{\sigma(x_i)},
\sigma_e^{-1}\xi_{\sigma(e)}
\rangle\\
&=\langle \sigma_e R(G,\sigma(p))\sigma_v^{-1}\xi_{\sigma(x_i)},
\xi_{\sigma(e)} \rangle
\end{align*}
and so $R(G,p)$ and $\sigma_eR(G,\sigma(p))\sigma_v^{-1}$ have the
same entry in the $\sigma(e)$ row and $\sigma(x_i)$ column.
Similarly, all entries agree.
\end{proof}

\subsection{Framework symmetries.}
Let $(G, p)$ be a framework in $\bR^d$ which is \textit{proper}
in the sense that the framework points are all distinct. Then a \textit{framework
symmetry} is a graph automorphism $\sigma$ of $G$ with the
additional property
\[
|p_{\sigma(i)}- p_{\sigma(j)}|= |p_{i}- p_{j}|
\]
for all edges $(v_i, v_j)$, where $|p_i - p_j|$ denotes Euclidean
distance.
%The \textit{framework symmetry group} of $(G,p)$ is the
%group of all such symmetries.
Note that such a symmetry  may just act locally. The framework of Figure 1
for example has  such a symmetry which exchanges $p_1$ and $p_3$.
We shall mainly be concerned
with the stricter \textit{global}
symmetries of frameworks that are
determined by isometric maps  of the ambient Euclidean
space. Thus we formally define a \textit{spatial symmetry}
 of $(G,p)$ as  a framework symmetry which is effected by an isometric map $T: R^d \to
\bR^d$ in the sense that
\[
\sigma(p) = Tp :=(Tp_1, \dots ,Tp_n)
\]
and we let $\G$ denote the \textit{spatial symmetry
group} of all such symmetries.
In the final section however we shall relax this and consider spatial symmetries in
subframeworks, and also latent spatial symmetries that appear after a partitioning.

The framework of Figure 1 has two evident mirror symmetries
which are spatial symmetries and $\G$ is isomorphic to the four group
$C_2 \times C_2$.

Recall that an isometric map $T$ admits a factorisation as
a product $T = T_1ST_2$, where $T_1, T_2$ are translations and $S$
is a linear isometry. The linearity of the entries in the rigidity
matrix ensures that $R(G,p) = R(G, Xp)$ if $X$ is a translation,
and it follows that $R(G,{T}p)$ = $R(G,{S}p)$. Consider $S$ also
in terms of the $d \times d$ real orthogonal matrix which effects
the transformation $p_i \to Sp_i$ by \textit{right} matrix
multiplication. In fact this matrix is $S^{-1}$ (where $S$ denotes
also the matrix that effects the transformation $S$). For example,
in case $d=2$, writing $(x_i',y_i')$ for the image $Sp_i$ of $p_i$
under $S$, we have
\[
\begin{bmatrix} x_i' & y_i'
\end{bmatrix} = \begin{bmatrix}x_i & y_i \end{bmatrix}S^{-1}.
\]
It follows from linearity that
\[
\begin{bmatrix}(x_i' - x_j')&(y_i' -y_j')\end{bmatrix} =
\begin{bmatrix}(x_i - x_j)& (y_i
-y_j)\end{bmatrix}S^{-1},
\]
and so
\begin{equation}\label{e:framesymmetry}
R(G,\sigma(p)) = R(G, Tp) = R(G, Sp) = R(G, p)\tilde{S}^{-1}
\end{equation} where $\tilde{S}= S \oplus \dots \oplus S$ is the
block diagonal matrix transformation of $\H_v$.

We now have all the ingredients for the proof of the individual symmetry equation
of part (i) of Theorem \ref{t:symm1}. For the general formula of part (ii)
we now specify five representations of the spatial symmetry group $\G$.

Write $\rho_e : \G \to \L(\H_e)$ for the permutation
representation of  $\G$ where
$\rho_e(g)$ is the transformation and the matrix which
is associated as
above with the spatial symmetry $g$.
Define $\rho_v : \G \to \L(\H_v)$ similarly. Let $\rho_{sp} :
\G \to \L(\bR^d)$ be the orthogonal group representation
(one often identifies $\G$ with
its image under this map) and let $\tilde{\rho}_{sp} = \rho_{sp}
\oplus \dots \oplus \rho_{sp}$ ($n$ times) be the associated block
diagonal representation of $\G$ on $\H_v$. Finally, note that
$\tilde{\rho}_{sp}$ and $\rho_v$ commute, that is,
\[
\tilde{\rho}_{sp}(g_1)\rho_v(g_2) = \rho_v(g_2)\tilde{\rho}_{sp}(g_1)
\]
for all $g_1, g_2$. Thus the product representation, denoted
$\hat{\rho}_v$, is well-defined. Indeed, these representations
can be viewed as representations
in different factors of the natural tensor product identification
$\H_v = \bR^n \otimes \bR^d$ and $\hat{\rho}_v= \rho_{n} \otimes
\rho_{sp}$, where $ \rho_n$ is the (multiplicity one)
representation for the vertices, so that $\rho_v = \rho_{n}\otimes
Id_d$, and $\tilde{\rho}_{sp}=Id_n \otimes \rho_{sp},$ where  $Id_n$ denotes the identity
representation of multiplicity $n$.

The next theorem provides  symmetry equations for the rigidity
matrix. For an alternative somewhat more sophisticated
derivation one may employ the chain rule for the
derivative of composite multi-variable functions and we do this in
Section 4 in a more abstract setting.

% and  follows from Lemma \ref{l:symm} and the discussion
%above.

\begin{thm}\label{t:symm1}
Let $(G,p)$ be a framework in $\bR^d$ with graph
$G=(V,E)$.

(i) If $T$ is a spatial symmetry for the framework
$(G,p)$ with associated graph symmetry $\sigma: V\to V$ and linear transformation
matrices $\sigma_v$ and $\sigma_e$
then
\[
R(G,p) = \sigma_e^{-1}R(G,p)\sigma_v\tilde{S}
\]
where $S$ is the linear isometry factor of $T$
and $\tilde{S}=S\oplus \dots \oplus S$ is the induced operator on $\H_v$.

(ii) Let $\G$ be the spatial symmetry group of the framework
$(G,p)$ with representation $\hat{\rho}_v = \rho_n \otimes
\rho_{sp}$ on $\H_v$ and representation $\rho_e$ on $\H_e$.
Then, for all $g \in \G$,
\[
R(G,p) = \rho_e(g^{-1})R(G,p)\hat{\rho}_v(g).
\]
\end{thm}

\begin{proof}
We may combine the equations \ref{e:graphsymmetry} and \ref{e:framesymmetry}
to obtain
\[
\sigma_e^{-1} R(G,p) \sigma_v = R(G,p)\tilde{S}^{-1},
\]
from which (i) follows. Now (ii) follows from (i) and the definition of the
representations $\rho_e$ and $\hat{\rho}_v(g)$.
\end{proof}

We note some immediate consequences for rigidity and isostaticity.

%\begin{definition}
%A  framework $(G,p)$ in the plane is infinitesimally rigid if
%the rank of $R(G,p)$ is  $2|V|-3$ and is isostatic if it is
%infinitesimally rigid and $|E| = 2|V|-3$.
%\end{definition}

The analysis above applies also to what one might call
\textit{grounded} or \textit{supported} frameworks $(G, p^*)$ in
which certain vertices are fixed absolutely. The relevant
symmetries in this case permute these special points. Such
examples can be found in the original three-point-supported
symmetric two-dimensional structures in  Kangwai and Guest\cite{kan-gue} and
Fowler and Guest\cite{fow-gue}.
\begin{center}
\begin{figure}[h]
\centering
\includegraphics[width=8cm]{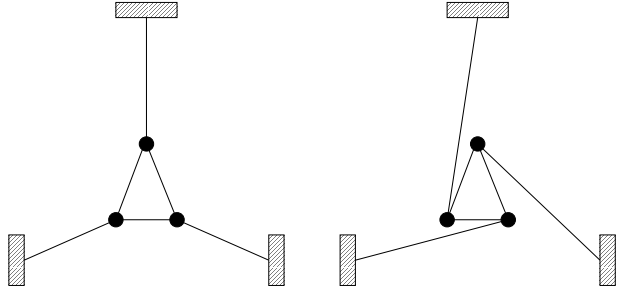}
\caption{The singular Jacobian for the first framework is a
consequence of reflection symmetry.}
\end{figure}
\end{center}
The context is simpler since  spatial flexes are absent and
isostaticity of the suspended framework corresponds to the
invertibility of the Jacobian $J(G,p^*)$ for the equation system
for the free points. The argument for Theorem 2.2 (ii) applies and
we obtain
\[
J(G,p^*) = \rho_e(g^{-1})J(G,p^*)\hat{\rho}_v(g),
\]
which is valid for elements $g$ of the spatial symmetry group $\G$, where
$\rho_v$ is the representation for free vertices. In particular if
$(G,p^*)$ is isostatic then
\[
\rho_e(g) = J(G,p^*)\hat{\rho}_v(g)J(G,p^*)^{-1}
\]
and so we obtain  the following equalities
of traces (also called characters); for each spatial symmetry group
element $g$,
\begin{align*}
\trace(\rho_e(g)) &= \trace(\hat{\rho}_v(g))\\
                  &=\trace(\rho_n(g)\otimes \rho_{sp}(g))\\
                  &=\trace(\rho_n(g))\trace(\rho_{sp}(g)).
\end{align*}
For the identity symmetry element one obtains the simple counting
condition $e' = 2v'$, where $e'$ is the number of bars and $v'$ is
the number of free joints. If a reflection symmetry $g=\sigma$
exists then since $\trace(\rho_{sp}(\sigma))=0$ one obtains
$\trace(\rho_e(\sigma))=0$ which is to say that there can be no
edges that are left fixed by the reflection.

As an illustration, consider the bilaterally symmetric frameworks of Figure 2.
From the above it follows that
$(G,p)$ is not isostatic if there is a reflection symmetry of the
framework which leaves invariant at least one edge. In this manner the symmetry equation serves as a device for recognising singular systems which is somewhat
simpler than the full Fowler-Guest equation.

\section{Flexes, Stresses and the Fowler-Guest Formula}

Let
$(G,p)$ be a proper framework in $\bR^d$, that is,
one with distinct framework points, and let ${\H}_{fl} = \ker R(G,p)$
and let ${\H}_{st} = \ker R(G,p)^*$ denote the kernel (nullspace)
of the adjoint (conjugate transpose) matrix. The notation reflects the fact
that the vectors of $\H_{fl}$
%(resp. $\H_{st}$) %
can be interpreted as infinitesimal \textit{flexes}
%(resp.internal stresses)
of the framework and that the vectors of $\H_{st}$
represent \textit{self stresses} of the framework, as we have indicated above.
In fact the infinitesimal flexes are
the vectors in the kernel of the derivative of
the nonlinear mapping from framework points coordinates to framework
edge lengths. This derivative, as we have noted is twice the rigidity matrix.

The symmetry equation shows immediately the key fact that for all $g \in {\G}$,
\[
\hat{\rho}_{v}(g){\H}_{fl} = {\H}_{fl},~~~~~~~
{\rho}_{e}(g){\H}_{st} = {\H}_{st}.
\]
That is, that these spaces are \textit{invariant subspaces} for the representations.
Thus with respect to the orthogonal decompositions ${\H}_{v}=
{\H}_{v'} \oplus {\H}_{fl}, \\
{\H}_{e}= {\H}_{e'} \oplus
{\H}_{st}$ the matrix $R$ takes the block form
$$R=\begin{bmatrix}
R'&0\\0& 0\end{bmatrix}
$$
where $R^{'}$ has trivial kernel and maps $\H_{v'}$ onto
${\H}_{e'}$. Certainly ${\H}_{fl}$ is nonzero since it contains
the space, $\H_{rig}$ say, corresponding to ambient rigid body
motion. In the case $d=2$ one may take as a basis for  $\H_{rig}$
the vectors
\[
u_{x}= (1, 0, 1, 0, \dots, 1, 0), u_{y}= (0, 1, 0, 1, \dots, 0,
1),
\]
(which are associated with infinitesimal translation), together
with the vector $u_{xy}$ (associated an infinitesimal rotation
about the origin) given by
\[
u_{xy} = (-y_1,x_1,-y_2,x_2, \dots ,-y_n,x_n),
\]
where $(x_i, y_i)$ are the coordinates of the framework
points $p_i$.
In fact for the associated three dimensional subspace
\[\H_{rig} =
{\H}_{x}\oplus \H_y \oplus {\H}_{xy}
\]
the subrepresentation $\rho_{rig}$ of $\hat{\rho}_v$ (obtained by
restriction to $\H_{rig}$) decomposes as $3$ copies of the trivial
one-dimensional representation.

Finally, define  $\H_{mech}$ as the complementary space  to
$\H_{rig}$ in $\H_{fl}$, so that $\H_{fl} = \H_{mech} \oplus
\H_{rig}$. The notation reflects the fact that this subspace may
be viewed as the space for non-trivial infinitesimal motions
(mechanisms) of the framework.

With these Euclidean space decompositions, which are all in terms
of invariant subspaces for $\hat{\rho}_v$), we have the associated
decompositions
\begin{equation}\label{e:rhovdecomp}
\hat{\rho}_{v}= \rho_{v^{'}} \oplus \rho_{fl} =  \rho_{v^{'}}
\oplus \rho_{mech} \oplus \rho_{rig}.
\end{equation}
For the other representation $\rho_e$ we have the two-fold
decomposition
\begin{equation}\label{e:rhoedecomp}
\rho_{e}= \rho_{e^{'}} \oplus \rho_{st}
\end{equation}
associated with the orthogonal decomposition ${\H}_{e}=
{\H}_{e^{'}} \oplus {\H}_{st}$.

We can now give a complete proof of a general form of the
Fowler-Guest formula \ref{e:FGequation}. In brief, the formula follows immediately
from the similarity (and unitary equivalence) of the "residual"
representations $\rho_{v'}$ and $\rho_{e'}$ following the removal
of subrepresentations corresponding to flexes (both ambient and
nontrivial) and to stresses, and this similarity follows from the
symmetry equation for the rigidity matrix.

Write $[\rho_{x}]$ to denote the character list  of a
representation $\rho_x$. Explicitly, this is the list
$(\trace(\rho_x(g_1)), \dots , \trace(\rho_x(g_N)))$ for some
enumeration of the elements (or the generators) of $\G$.

\begin{thm}\label{t:FG}
Let $(G,p)$ be a bar-joint framework in $\bR^d$, with $n$ distinct joints
and $m$ bars, and with spatial symmetry group $\G$ with
orthogonal representation $\rho_{sp}$ in $\bR^d$. Let $\rho_n$,
$\rho_e$  be the joint and bar (permutation) representations of $\G$ on
$\bR^n$ and $\bR^m$ respectively, let $\rho_{rig}$ be the
subrepresentation of $\rho_{n}\otimes \rho_{sp}$ for the space of
trivial infinitesimal flexes, let $\rho_{mech}$ be the
subrepresentation for nontrivial flexes, and let $\rho_{st}$ be
the subrepresentation of $\rho_e$ for the space of internal
stresses. Then
\begin{equation}\label{e:fowguesttracelist}
[\rho_{{mech}}]- [\rho_{st}]= [{\rho}_{n}]\cdot [\rho_{sp}] -
[\rho_{e}]- [\rho_{rig}]
\end{equation}
where $[~]\cdot[~]$ denotes entry-wise product of characters.
\end{thm}

\begin{proof}
Recall that with respect to the orthogonal decomposition
$${\H}_{v}=
{\H}_{v}^{'} \oplus {\H}_{fl},
{\H}_{e}= {\H}_{e}^{'} \oplus
{\H}_{st}
$$
the rigidity matrix $R$ takes the block form
$$R=\begin{bmatrix}
R'&0\\0& 0\end{bmatrix}.
$$
The matrix $R^{'}$ is
a square nonsingular matrix which we view as a linear
transformation $R' : \H_v' \to \H_e'$. From the
symmetry equation we have the commutation
relations $R^{'}\rho_{v^{'}}(g)= \rho_{e^{'}}(g)R^{'}$
and so
\[
\trace(\rho_{v'}(g))=\trace((R')^{-1})\rho_{e'}(g)R'=\trace(\rho_{e'}(g)).
\]
Thus $\rho_{v^{'}}$ and $\rho_{e^{'}}$
have the same character list; $
[\rho_{v'}]=[\rho_{e'}]$.

We have \begin{align*}
[{\rho}_{n}]\cdot [\rho_{sp}]=[{\rho}_{n}\otimes \rho_{sp}]&= [\rho_{v'}]+[\rho_{mech}]+[\rho_{rig}]
\end{align*}
and
\begin{align*}
[\rho_e] &= [\rho_e']+[\rho_{st}]
\end{align*}
%From the unitary equivalence induced by $R'$
and so from $[\rho_{v'}]=[\rho_{e'}]$ we obtain equation
\ref{e:fowguesttracelist}.
\end{proof}

We note that one can also make explicit an orthogonal equivalence between the residual representations $\rho_{v'}, \rho_{e'}$
through the  isometric part $U$  of the polar decomposition
$R' = U|R'|$ as this  operator
also intertwines the representations;
 $U\rho_{v^{'}}(g)= \rho_{e^{'}}(g)U$, for $g \in \G$.

The right hand side of the Fowler-Guest formula is readily computable in
terms of the number of elements fixed by a framework symmetry.
Thus one may quickly obtain necessary counting conditions for such elements
when the number of independent stresses and mechanisms are specified.
Recall that a framework $(G,p)$ is isostatic if it is infinitesimally
rigid and is stress free.
Thus, in the case of planar isostatic framework
if there is a mirror symmetry  $\sigma$ evaluation of the formula
at  $\sigma$ gives
\[
0-0=\trace(\rho_n(\sigma))\trace(\rho_{sp}(\sigma)) -\trace(\rho_e(\sigma))-\trace(\rho_{rig}(\sigma).
\]
Since $\trace(\rho_{sp}(\sigma)) =0$ and $\trace(\rho_{rig}(\sigma))=1$
we obtain
\[
0=0-b_\sigma+1
\]
where here we follow Connelly et al\cite{con-et-al} and write
$j_\sigma$ and $b_\sigma$ for the number of framework points
(joints) and framework edges (bars) that are not displaced by
$\sigma$.

Let us consider the framework of Figure 1 once more.
Adding a cross edge between one of $p_1,p_2,p_3,p_4$ and one of
$p_5,p_6,p_7,p_8$ will create  realisations of a Laman graph.
Addition of $[p_2,p_6]$ removes both mirror symmetries, so the count condition
above is irrelevant and indeed the framework
is isostatic. Also note that addition of $[p_2,p_6]$ is consistent with the necessary condition for the one remaining mirror symmetry. On the other hand addition of an edge on the $x$-axis violates the count $b_\sigma =1$, for both mirror
symmetries, and for this reason the resulting framework is not isostatic.

We can obtain also the following corollary which is indicative of the results
obtained in Connelly et al\cite{con-et-al}.

\begin{cor}\label{countingCor}
Let $(G,p)$ be an isostatic framework in $\bR^3$ which does not
lie in a hyperplane and which has a proper spatial symmetry
$\sigma$. Then the following equations hold.

(i) If $\sigma$ is a half turn then $0 = -j_\sigma-b_\sigma+2$.

(ii) If $\sigma$ is a reflection then $0=j_\sigma-b_\sigma$.

(iii) If $\sigma$ is an inversion then $0 = -3j_\sigma-b_\sigma$.
\end{cor}

\begin{proof}
(i) In this case $\trace(\rho_n(\sigma))= j_\sigma$ since
$\rho_n(\sigma)$ is a permutation matrix with a nonzero diagonal
entry if and only if the corresponding vertex is fixed by
$\sigma$. Also $\trace(\rho_{sp}(\sigma))=-1,$ since
$\rho_{sp}(\sigma)$ is equivalent to a diagonal matrix with
entries $-1,-1,1$, and $\rho_{rig}(\sigma)= -2$ since in the three
dimensional subspace for infinitesimal translation flexes
$\rho_{rig}(\sigma)$ is diagonal with entries $-1,-1,1$, and in the
three dimensional subspace for infinitesimal rotation flexes
$\rho_{rig}(\sigma)$ is similarly diagonal with entries $-1,-1,1$.

From these observations and the previous character formula,
evaluated at $\sigma$, statement (i)  follows. The formulae of (ii) and
(iii) are similarly verified; in case (ii),
$\trace(\rho_{sp}(\sigma)) =1, \trace(\rho_{rig}(\sigma))=0$ and in
case (iii), $\trace(\rho_{sp}(\sigma))=-3$ and
$\trace(\rho_{rig}(\sigma))=0$.
\end{proof}

\section{Higher Order Frameworks and Symmetry}

We now show how the approach above adapts readily to higher
dimensional frameworks such as point-line frameworks in $\bR^2$,
body-bar frameworks in $\bR^3$, and even infinite frameworks.

\subsection{Character formulae for point-line frameworks.}
 Consider, in $\bR^2$, a set $\P$ of points and a
set $\L$ of straight lines,
\[
\P =\{p_1, \dots , p_n\}, \quad\L = \{L_1, \dots , L_r\}.
\]
%In the first instance we view the lines as directed.
Considering only certain pairs from $\P \cup \L$ we can compute
generalised distances involving the lines, namely point-line
distances, being the usual nonnegative distance, and line-line
angles, taking values in $[0,\pi /2]$. The chosen pairs determine
edges $e\in E$ in an abstract graph whose vertex set $V$ is
partitioned $V = V_p \cup V_l$  and whose edge set is similarly
partitioned, $E=E_{pp}\cup E_{pl}\cup E_{ll}$. The abstract
partitioned graph $G$ and the pair $\P, \L$ give rise to a
distance labeled graph. This is the pair $(G, d)$ where $d$ is a
map from $E$ to the set of distances; $d(e) = d(p_i,p_j),$ for $
e=(i,j)\in E_{pp}, d(e) = d(p_i,L_j), $ for $ e \in E_{pl} $ and $
d(e)= d(L_i,L_j)$ for $ e=(i,j) \in E_{ll}$.

It is of interest to understand the inverse problem, that is, the
nature of solutions of the constraint equations determined by an
abstract  distance labeled partitioned graph. These equations are
in the coordinate  variables for the points and lines.
The
points are coordinatised as usual, with variables $x_i, y_i$ for
the framework point $p_i$. We may assume by
translating that the lines $L_j$ do not pass through
the origin and so may be parameterized by their closest points
$(x_j', y_j')$ to the origin.
Writing $\ul{x}$ for the set
of all variables, this system can be indicated as the equation set
\[
f_e(\ul{x}) = d(e), \quad e\in E.
\]

Let $(G, \P, \L)$ be a point-line framework as above. Define
$\H_e$ and $\H_v$  as before but with the natural additional
structure:
\[
\H_v = \H_p \oplus \H_l
\]
and, according to edge type,
\[
\H_e = \H_{pp}\oplus \H_{pl}\oplus \H_{ll}.
\]

 Also,
\[
{\H}_{p}= \sum_{k=1}^{n}\oplus (\mathbb{R}_{x_{k}} \oplus
\mathbb{R}_{y_{k}}),\quad {\H}_{l}= \sum_{k={n+1}}^{n+r}\oplus
(\mathbb{R}_{x_{k}'} \oplus \mathbb{R}_{y_{k}'}).
\]

We define the rigidity matrix for a line-plane framework, or a
dimensioned abstract graph, simply  as the Jacobian of the
distance constraint equation system.
The Jacobian has a $3 \times 2$ block structure implied by the
vector space decompositions and takes the form,
$$R(G,\P,\L)=\begin{bmatrix}
R(G,\P)&0\\0& R(G, \L)\\
R_1&R_2
\end{bmatrix},
$$
and the representations $\rho_e$ and $\hat{\rho}_v$ have a
corresponding three-fold and two-fold diagonal block structure,
respectively.

As before we have a spatial symmetry group $\G$ for the point-line framework
$(G,\P,\L)$. For simplicity we assume that the framework contains lines and points, that $0 \in \bR^2$ is the centre of symmetry and that there are no lines through the origin.
As before we have five representations  :
\[ \rho_e,
~~\rho_v, ~~\rho_{sp}, ~~\tilde{\rho}_{sp} ~~\mbox{ and }~~
\hat{\rho}_v = \rho_n \otimes \rho_{sp}.
\]
Note in particular that the coordinates for the lines are analogous to
the coordinates for points in that for a point-line framework symmetry
$g$, given by a linear isometric transformation $T$ of $\bR^2$, the
coordinates for the transformed line $T(L_j)$ are $T(x_j',y_j')$.

Define $\H_{st} = coker R(G,\P,\L)$ and let $ker
R(G,\P,\L) = \H_{rig} \oplus \H_{mech}$ where $\H_{rig}$ is the
three dimensional space of infinitesimally rigid flexes. Thus the space
$\H_{mech}$ is defined as the (possibly zero) orthogonal
complement of $\H_{rig}$ in $ker R(G,\P,\L)$.
The point-line framework is said to be \textit{isostatic} if it is infinitesimally
rigid, that is, if if $\H_{mech}=\{0\}$, and also that it is stress free
in the sense that $\H_{st}=\{0\}$.

\begin{thm}\label{t:pointline}Let
$(G,\P,\L)$ be a point-line framework as above with spatial symmetry group
$\G$. Then
\[
R(G,\P,\L) = \rho_e(g^{-1})R(G,\P,\L)\hat{\rho}_v(g), \quad g\in
\G,
\]
and, as an equality of character lists,
\[
[\rho_{{mech}}]- [\rho_{st}]= [{\rho}_{n}]\cdot [\rho_{sp}] -
[\rho_{e}]- [\rho_{rig}].
\]
In particular if the framework is isostatic and has a proper
reflection symmetry, with graph automorphism $\sigma \neq id$,
then
$$ b_{pp} + b_{ll} +b_{pl} -
1 =0 $$ where  $b_{pp}$, $b_{ll}$ and $b_{pl}$ are the number of
point-point edges, line-line edges and point-line edges which are
unchanged by the reflection.
\end{thm}

\begin{proof}
In the next subsection we obtain a  general symmetry formula and
the stated formula is a special case of this. The character list
formula is proven in exactly the same manner as in the proof of
Theorem \ref{t:FG}
\end{proof}

As we have noted in the introduction,
this theorem can be useful for predicting the singularity of an
equation system underlying a CAD diagram.
%For a simple
%illustration of this consider the triangular point-line drawing
%ramework of Figure 3. The abstract graph has six vertices and
%nine constraints and no subgraphs are overconstrained. However
%there is one reflection symmetry, for which $b_{pp}=1, b_{ll}=1$
%and it follows that the Jacobian of the equation system is
%singular.

\subsection{Higher order frameworks}
We now derive symmetry equations for the
rigidity matrix of quite general distance constrained systems using a more direct
proof using the Jacobean derivative of the generalised
edge map.
%Again one can obtain
%character formulae and  counting criteria for isostatic
%systems.
A simple example of the abstract formulation below is the case of finite systems of points
and (unoriented) planes in $\bR^3$, with constraints of Euclidean
distance between points, and points and planes, and with angular
constraints between planes. Planes may be coordinatised by the
three coordinates of the point closest to the origin and so play a
role similar to points.
%IN this case, with $j_\sigma $ (resp. $b_\sigma$)
%counting the total number of points and planes (resp. constraint
%edges) left undisplaced by $\sigma$ one can recover the necessary
%conditions of Corollary 3.2.}

Let $(G,E)$ be a finite, connected, undirected graph and let
$V=V_1\cup \dots \cup V_n$ be a partition in which $V_i=\{v_{i,k}
:1 \leq k \leq \nu_i\}$ is a set of vertices which label a set
$\P_i =\{p_{i,k}:1 \leq k \leq \nu_i\}$ of geometric objects of
the same kind. Formally, each \textit{object} of the $i^{th}$
kind, $p_{i,k} \subseteq \bR^d$, is a real manifold, or, more
generally, a real semi-algebraic set, which is determined by a
specification $\underline{x}_i= (x_{i,1},\dots ,x_{i,t_i})$ of
$t_i$ parameters. For example, a straight line in three dimensions
requires four variables.
%The equation system for a particular
%object follows a choice of coordinatisation  and only plays a role
%in the determination of the constraint function for the parameters
%of the objects.

For a pair of specified objects $(p, q)$, either of the same or
differing type, a generalised distance equation is given which has
the form $f(p,q)=d$ where $d$ is real and $f$ is a function in the
parameters for $p, q$.
%This (general) distance
%equation also gives rise to a a constraint equation for
%unspecified objects of these types, namely $f(p,q)=d$ with $d$
%viewed as input rather than output variable.
We  say that the constraint is a \textit{Euclidean constraint} if
for all isometries of $\bR^d$ and all objects $p, q$ of the
appropriate type, we have $f(Tp,Tq) = f(p,q)$.
%(ADD COMMENTS on
%quadratic pols with no constant term??).

\begin{definition}
A {Euclidean framework} is a pair $(G, \P)$ together with a family
of distance functions $f_e(p,q), e \in E$ where

(i) $G=(V,E)$ is a graph with partitioned vertex set $V$ labeling a
set $\P$ of specified objects, with objects of the same kind in
each  partition set, and

(ii) the distance functions $f_e(p,q)$ are Euclidean invariant and
depend only on the type of the objects $p, q$.
\end{definition}

To consider the rigidity or flexibility of a particular Euclidean
framework $(G, \P)$ we consider  the \textit{framework equation
system}, which, by definition, is the constraint system
\[
f_e(\ul{x}_{i,k},\ul{x}_{j,l}) = d_e, \quad e=(v_{i,k},v_{j,l})
\in E,
\]
A proper Euclidean framework $(G, \P)$ is one for which the
objects do not all lie in a hyperplane. We say that  a
framework of this type is  \textit{infinitesimally rigid} if the Jacobean
$J(G,\P)$ of the constraint system has rank equal to $N-d(d+1)/2$ where
\[
N=\nu_1t_1+\dots +\nu_nt_n
\] is the total number of variables. Also we say that $(G, \P)$ is
\textit{isostatic} if in addition the rank is equal to $|E|$.

Let $(G, \P)$ be a Euclidean framework  with geometric objects
$p_1, \dots ,p_n$. Following the terminology for frameworks we
define the constraint function, or edge  map, of $(G,\P)$ to be
the nonlinear function $f : \bR^N \to \bR^m $ with
\[
f(\ul{x})=(f_{e_1}(\ul{x}),\dots ,f_{e_m}(\ul{x})),
\]
where $m=|E|$.
Here the $i^{th}$ constraint function for the edge $e_i$ depends
on the variables $\ul{x}_k, \ul{x}_l$ for the objects $p_k, p_l$
associated with $e_i$.

We have
\[
\H_v = \sum_{k=1}^{n}\sum_{i=1}^{\nu_k} \oplus \bR^{t_i},
\]
as the vector space for coordinate variables
and if $T$ is an isometric transformation of $\bR^d$ then there is
an associated block diagonal transformation
\[
\tilde{T}= \sum_{k=1}^{n}\sum_{i=1}^{\nu_k} \oplus T_k,
\]
where each $T_k$ is the parameter transformation induced by $T$.
In particular, if $\sigma$ is a spatial symmetry of $(G, \P)$ which
additionally is induced by a spatial isometry $T$ then we call
$\tilde{T}$ the local symmetry transformation for $\sigma$.

Similarly we have the edge space $\H_e$ on which the spatial symmetries $g$ act
a permutation transformations.

We now obtain the symmetry equation for the rigidity matrix of a
Euclidean framework, defined here as Jacobean derivative
$D(f)(\ul{x})$ of the constraint map evaluated at the framework
coordinates to yield the matrix $J(G, \P)$.

\begin{thm}\label{t:euclidean} Let $(G, \P)$ be a Euclidean framework, with
generalised distance equations $f_e(p,q)= d_e, e\in E$, where $p,
q$ denote the parameters of the two geometric elements constrained
by distance $d_e$, let $f:\bR^N \to \bR^m$ be the generalised
constraint  map and let $(\sigma, T)$ be a spatial symmetry of
$(G,\P)$. Then the rigidity matrix $J(G, \P)$   satisfies the
symmetry equation
\[
J(G, \P)=\sigma_e^{-1}J(G, \P)\sigma_v\tilde{T} =\sigma_e^{-1}J(G,
\P)\tilde{T}\sigma_v.
\]
where $\sigma_v$ and $\sigma_e$ are the induced permutation
transformations of the vertex space $\H_v$ and the edge space
$\H_e$ and where $\tilde{T}$ is the local symmetry transformation
for $\sigma$.
\end{thm}

\begin{proof}
Let $\sigma$ and $T$ be as above. Then from the graph symmetry
$\sigma$ it follows, as in Lemma \ref{l:symm}, that evaluating the Jacobian
at $\sigma(x)$ gives the same matrix as corresponding row and
column operations on the Jacobian, that is,
\[
%Df(\tilde{T}\ul{x})=
Df(\sigma(\ul{x})) = \sigma_e^{-1}Df(\ul{x})\sigma_v.
\]
On the other hand, by Euclidian invariance $f(\tilde{T}\ul{x}) =
f(\ul{x})$ for all values of the variables, and so by the chain
rule,
\[
(Df)(\tilde{T}\ul{x})\tilde{T} = Df(\ul{x}).
\]
However, we have $\sigma(\ul{x}) = \tilde{T}\ul{x}$  for the given
framework coordinates and putting these fact together yields in
this case
\[
Df(\ul{x})\tilde{T}^{-1}= Df(\tilde{T}\ul{x})= Df(\sigma(\ul{x}))
= \sigma_e^{-1}Df(\ul{x})\sigma_v,
\]
as required.
\end{proof}

%\begin{rem}{\rm One can again use the general rigidity matrix symmetry formula to obtain
%character formulae and hence counting criteria for isostatic
%systems. A simple example is the case of finite systems of points
%and (unoriented) planes in $\bR^3$, with constraints of Euclidean
%distance between points, and points and planes, and with angular
%constraints between planes. Planes may be coordinatised by the
%three coordinates of the point closest to the origin and so play a
%role similar to points. With $j_\sigma $ (resp. $b_\sigma$)
%counting the total number of points and planes (resp. constraint
%edges) left undisplaced by $\sigma$ one obtains the necessary
%conditions of Corollary 3.2.}
%\end{rem}

\subsection{Pin-jointed body frameworks}
 We now consider a generalisation of  bar-joint frameworks by
allowing the edges to be general rigid bodies which may then have
more than 2 vertices. Informally this looks like a set of rigid
bodies which are held together by a set of pins or hinges, each of
which passes through two or more bodies. Note that bar-joint and
body-bar frameworks are both special cases of pin-jointed body
frameworks. The discussion below is self-contained. For other
information on body bar frameworks see Tay and Whiteley
\cite{tay-whi} and Jackson and Jordan\cite{jac-jor}.
We limit attention to pin-jointed body frameworks in $\bR^2$.

\begin{definition}
A pin-jointed body framework is a pair $(\S,p)$ where $p =
\{p_i\}$ is a set of points in $\bR^2$ and $\S = \{S_e\}$  a collection of
subsets of the points such that:

(i)  every point is in at least two sets,

(ii) every set contains at least two points.
\end{definition}

We also shorten the term to "body framework" and we call the sets
$S_e$ "bodies". The labelling notation here reflects the special
case of edges and we occasionally denote a set  $S_{e_i}$ simply by $e_i$.
Every body framework defines a bipartite graph $G=G(\S)$ in which
the points are the vertices of one partition and the bodies are
the vertices of the other partition. The edges of $G$ represent
the occurrence of a point in a body. Conversely a bipartite graph
with minimum vertex degree  greater than one defines a
body framework.

A \textit{flex} (or infinitesimal flex, or  infinitesimal motion)
of a body framework is an assignment of velocities $u_i$ to the
points $p_i$ and an assignment of infinitesimal motions $(v_e,
a_e)$ to the bodies such that for each body the velocities of its
points are compatible with the rigid motion $(v_e, a_e)$ of the
body. Here $v_e\in \bR^2$ is the velocity of the centroid of the
body $e$ and $a_e\in \bR$ is its angular velocity, and the
centroid is defined as $p_e = \frac{1}{| S_e |}\sum_{p_i \in
S_e}p_i$. The compatibility condition is the equation
\[
u_i = v_e + a_e(p_i-p_e)^{\pi/2},
\]
where  $v^{\pi/2}$ denotes the rotated vector $(-y,x)$ when
$v=(x,y)$. Thus there are two linear equations for every
occurrence of a point in a body, that is, for every edge of the bipartite graph. With the coordinate notation $u_i
= (u_i(x),u_i(y))$ they take the form.
\[
u_i(x) - v_e(x) + a_e(p_i(y)-p_e(y))=0,
\]
\[
\quad u_i(y) - v_e(y) - a_e(p_i(x)-p_e(x))=0.
\]
 Suppose now that there are $n$ points, $e$ bodies and $c$
point-body occurrences, that is, $n+e$ vertices and $c$ edges in
$G(\S)$. We define a $(2n+3e)$ by $2c$ rigidity matrix $R= R(\S, p)$
as follows.

(i) $R$ has $2$ columns for each point and $3$ columns for each
body.

(ii) $R$ has 2 rows for each point-body occurrence.

(iii) The $2$ by $5$ submatrix for (i) and (ii) with appropriately
labeled columns, takes the form
%The two rows for the occurrence of point $p_i$ in body $e$
% contain the
%following entries. In the two columns for the point $p_i$, both
%rows have an entry 1 and an entry zero. In the first two of the
%three columns for the body e both rows have entry -1. In the third
%of these three columns one row has entry $p_i(y)-p_e(y)$ and the
%other has entry $-(p(i(x)-pe(x))$. The $2$ by $5$ submatrix thus
%has the form
$$\begin{bmatrix}
u_i(x)   &u_i(y)     &v_e(x)    &v_e(y)     &a_e   \\
\hline \\
1&0&-1&0 &-(p_i(y)-p_e(y)) \\
0&1&0&-1&(p_i(x)-p_e(x))
\end{bmatrix}
$$

A body framework is \textit{infinitesimally rigid} if it has no
non-trivial flexes. As usual there is a three-dimensional space of
trivial flexes and so infinitesimal rigidity corresponds to there
being no other nonzero solutions to the compatibility equations.
This is simply the condition ~~~dim(ker$R) =3$. We say that a body
framework is \textit{isostatic} if $2c = 2n+3e-3$ and rank $R =
2c$.

Consider now the  natural decompositions of the domain space
and the codomain space for the rigidity matrix regarded as a
linear transformation.

Let $p_1, \dots , p_r$ be the pin points of $(\S,p)$ and let
$e_1,\dots ,e_s$ be the bodies. Let
$
\H_{dom} = \H_{b}\oplus \H_{p},$ where
\[
\quad  \H_{b}=
\H_{body}\otimes \bR^3=
\sum_{i=1}^s
\oplus \bR^3,
\quad \H_p = \H_{pin}\otimes \bR^2 = \sum_{i=1}^r \oplus \bR^2,
\]
where the summands $\bR^2$ represent the spaces of displacement
velocities $u_i$ for $p_i$ and where the summands $\bR^3$ are the
spaces of body velocities $(v_e(x), v_e(y),a_e)$. Similarly, the codomain
space for $R$ has the form
\[
\H_{codom} = \H_{mem}\otimes \bR^2=\sum_{i=1}^N \oplus \bR^2,
\]
associated with the $N$  edges of the bipartite graph of $(\S,p)$,
that is, with the membership conditions $p_i \in e_j$.

Let $\G= \G(\S,p)$ be the group of isometries $T$ of $\bR^2$ that
are body-framework symmetries. Thus $Tp_i = p_{\pi(i)}$ for some
permutation $\pi$ of the pins, and $\pi$ respects bodies, that is,
the set $\pi(e_i)$ is equal to $e_{\tau(i)}$ for some permutation
$\tau$. In particular the pair $(\pi, \tau)$ gives an automorphism
of the abstract bipartite graph of the body framework.

Once again we consider various natural representations of $\G$.
First we have  $\rho_{b}=\rho_{body}\otimes Id_3$ and $ \rho_{p}=\rho_{pin}\otimes Id_2$, the (inflated) permutation
representations of the spatial symmetry group $\G$ on $\H_{b}$
 and $\H_{p}$ associated with $\pi$ and $\tau$ respectively.
As before, let $\rho_{sp}$ be the spatial representation of $\G$ as
orthogonal transformations of $\bR^2$, and let $\rho_{sp}^+$
be the representation $\rho_{sp}\oplus \Delta$ on $\bR^3$ where
$\Delta$ is the one dimensional determinant representation. We then have the natural  representation of $\G$ on
$\H_{dom}$ given by
\[
\rho_{dom} := \hat{\rho}_{b}\oplus \hat{\rho}_{p}:=
(\rho_{body}\otimes \rho_{sp}^+) \oplus (\rho_{pin} \otimes \rho_{sp})
\]
where $\rho_{body}$ and $\rho_{pin}$ are the basic permutation representations
for bodies and for pins.

Secondly, there is a representation $\rho_{codom} = \rho_{mem}\otimes \rho_{sp}$
of $\G$  associated with the permutation representation $\rho_{mem}$ for
the edges of the bipartitie graph

In view of the form of the $2$ by $5$ submatrices above direct
calculation gives the symmetry equations
\[
 R = \rho_{codom}(g^{-1})R{\rho}_{dom}(g), \mbox{  for } g \in \G.
\]
As before these equations give to the invariance of various subspaces under the
representations $\rho_{dom}$ and $\rho_{codom}$;
$\rho_{rig}$ is the
subrepresentation of $\rho_{dom}$ determined by restriction to the
subspace $\H_{rig}$ of trivial rigid body motion flexes,
%This representation is
%simply three copies of the trivial one dimensional representation.
$\rho_{mech}$ is determined by the restriction to
$\H_{mech}:=$ ker$R\ominus \H_{rig}$, and $\rho_{st}$ is the
restriction of $\rho_{codom}$ to the (internal stress) subspace
$\H_{st} := coker R$.

\begin{thm}\label{t:body}
Let  $(\S,p)$ be a body framework in $\bR^2$ with spatial symmetry group
$\G$. Then the representation character lists satisfy the equation
\[
[\rho_{{mech}}]- [\rho_{st}]= [{\rho}_{sp}^+]\cdot [\rho_{body}]
+[\rho_{sp}]\cdot[\rho_{pin}] - [\rho_{codom}]- [\rho_{rig}].
\]
\end{thm}

\begin{proof}
The restriction of $R$ to the subspace $\H_{dom}\ominus (\H_{mech}\oplus \H_{rig})$
gives a linear bijection to $\H_{codom}\ominus \H_{st}$ and so the associated ("residual") representations are equivalent. The formula now follows, as in the proof of Theorem \ref{t:FG}
\end{proof}

As a corollary we see that if the body framework is isostatic and
has a reflection symmetry $\sigma$ then
\[
0=n_{body}^\sigma -1,
\]
where $n_{body}^\sigma$ is the number of bodies left unmoved by
$\sigma$. Indeed this follows from evaluating the character list equation at $\sigma$,
for we then have
\[
 \trace({\rho}_{sp}^+(\sigma)) =1, \trace(\rho_{body}(\sigma
)) = n_{body}^\sigma,
\]
as well as
\[ \trace(\rho_{sp}(\sigma)) =0,
~~\trace(\rho_{codom}(\sigma))=0, \mbox{ and } \trace(\rho_{rig}(\sigma))=-1.
\]

\subsection{Symmetry equations for infinite frameworks.}
In Owen and Power\cite{owe-pow-2,owe-pow-kempe,owe-pow-crystal} we have indicated some perspectives for a mathematical theory of infinite bar-joint frameworks. Part of the
motivation for such a development also comes from materials analysis
(Donev and Torquato
\cite{don-tor}), the analysis of repetitive structures
(Guest and Hutchison
\cite{gue-hut})
and from applications in chemistry (Ceulemans et al\cite{ceu-et-al}
and crystallography (Borcea and Streinu\cite{bor-str}).
%(See also\cite{owe-pow-kempe} and\cite{owe-pow-cry}.
We now consider the rigidity matrix symmetry equations  in this setting.
In particular we give a Hilbert space variant of Theorem \ref{t:FG}
for a natural notion of square-summable isostaticity,
and we give a  Fowler-Guest formula for
periodic frameworks.
 Of course a novelty for infinite frameworks is that the spatial
symmetry group $\G$ can be infinite.

\subsubsection{Infinite frameworks.} Let $(G,p)$ be a countable (and nonfinite)  bar-joint framework
in $\bR^2$ associated with a countable connected graph $G$, where
the framework vector $p=(p_1,p_2,\dots )$ has framework points
$p_i$ in $\bR^2$ indexed as usual by the vertices of $G$.
The consideration of such infinite  frameworks of a general
character, without translation symmetries,  was begun in Owen and Power\cite{owe-pow-2}.
Here the divergence of various notions of rigidity
was indicated  as well as forms of rigidity allied to  operator
interpretations of the rigidity matrix. This latter theme is developed further in
Owen and Power\cite{owe-pow-crystal}.
In addition to tools from  operator theory it seems that general notions from functional analysis (such as uniform convergence, compactness, aperiodicity) will become of relevance to the analysis of infinite framework deformability.
For our present consideration we address only infinitesimal rigidity
rather than continuous rigidity  and so we need only restrict attention to the rigidity matrix and its interpretations as a linear transformation.

Define the
rigidity matrix $R(G,p)$  as in Section 2.1, with the rows
labeled by edges and the columns labeled by vertices (twice over,
for $x$ and $y$ coordinates).
Assume that each vertex has finite degree. This entails that each column
of the matrix has finitely many nonzero entries.
This rigidity matrix may be viewed
as a linear transformation $T$ from the direct product vector
space $\H_v=\Pi_V \bR^2$ to the vector space $\H_e=\Pi_E  \bR$. Here
the direct product notation $\Pi_E
\bR$ indicates  the set of \textit{all} real sequences indexed by the edges
of $G$, with the usual vector space structure.
The permutation representation $\hat{\rho_v}$ and ${\rho_e}$ are defined on the spaces
$\H_v$ and $\H_e$, respectively, as before.

\begin{thm}\label{t:infsymm}
Let $(G,p)$ be an infinite bar-joint framework in $\bR^d$
with rigidity matrix transformation $R(G,p): \H_v \to \H_e$.
Then
$$R(G,p) = \rho_e(g^{-1})R(G,p)\hat{\rho}_v(g),\quad  g \in \G.$$
\end{thm}

\begin{proof}
The sparse nature
of the matrix for $R(G,p)$ ensures that
the various infinite sums implied by matrix multiplication
are  sums over finitely many nonzero terms.
With this change only the proof follows that of Theorem \ref{t:symm1}.
\end{proof}

Once again, we may choose
three linearly independent vectors  in the kernel of $T$ to span the linear
subspace of rigid motion  flexes associated with a three-dimensional space
$\H_{rig}$ for translations
and rotations.

 It is also natural to consider
$R(G,p)$ as a linear transformation between other smaller
sequence spaces which are invariant for the representations,
and in this case the symmetry equations will hold as above.
For example, let $T_0$ be the restriction of
$R(G,p)$ to the vector space direct sum, $\H_0=\Sigma_V \oplus \bR^2$,
which consists of \textit{finite} linear combinations of the usual standard basis
vectors ($\xi_{x_i}$ and $\xi_{y_i}, i=1,2,\dots )$. These are the
"finitely supported vectors", that is, the sequences $u=(u_v)_{v\in V}$ in $\H_0$ which have all but finitely many entries equal to zero.
One may view the vector $u$ as an assignment  of velocity
vectors to a finite number of joints of the infinite framework and view $T_0$ and associated mathematical constructs as
modeling a very large system and its finitely acting disturbances.
Note that $T_0$ maps into $\Sigma_E \oplus \bR$, in view of the finiteness of vertex degrees. Also note that the
translation and rotation  flexes do not lie in the domain
of $T_0$.
% If there is a nonzero vector $u$ in
%the kernel of $R(G,p)_0$ this is necessarily a nontrivial
%infinitesimal flex, and
It is natural then to say that $(G,p)$ is \textit{finitely
infinitesimally rigid} if the kernel of $T_0$ is trivial.  The
regular square grid framework (with framework points $(i,j), i,j
\in \bZ$) has this property as do grid frameworks with more generic
vertex locations. Indeed it is enough to show that for
any finite large square grid there is no nonzero flex which assigns zero
velocities to the boundary joints.
%as does the Kagome framework of\cite{gue-hut}, and
In fact we say that this framework is  \textit{finitely isostatic}
since in this case there are also no nontrivial finitely supported
stresses (vectors in the cokernel).

One can also consider other less severe constraints on the domain
space, that is, on the allowable velocity vectors and flexes $u$,
such as   boundedness (each domain vector $u$ is a bounded
sequence), summability ($\sum_v |u_v| < \infty$), or square
summability ($\sum_v |u_v|^2 < \infty$).

Let us define a \textit{square-summably isostatic
framework} in $\bR^d$ as one for which

(i) the rigidity matrix $R(G,p)$ determines a bounded Hilbert
space operator $T(G,p)$ from the real Hilbert space
$\H_v^2:=\ell^2(V)\otimes \bR^d$ to the real Hilbert space $\H_e^2:=\ell^2(E)$,

(ii) the kernel and cokernel of $T(G,p)$ are the zero subspaces.

Once again, for the spatial symmetry group  we
have the representations $\hat{\rho_v} =\rho_v\otimes \rho_{sp}$,
on $\H_v^2$ and $\rho_e$ on $\H_e^2$.
The following proposition is  an
infinite framework generalisation of the unitary equivalence
noted in the finite case for the residual representations of $\G$.

\begin{prop}
Let $(G,p)$ be a square-summably isostatic framework in $\bR^d$.
Then
$\hat{\rho_v} $ and $\rho_e$ are unitarily equivalent representations and in particular
have the same irreducible components.
\end{prop}

\begin{proof}
We use a standard argument to show that the unitary part
of $T=T(G,p)$ implements the equivalence.

Since $(G,p)$ is square summably isostatic $T$ has a unique polar
decomposition of the form $T=U|T|$ with $U$  unitary. We have
$\rho_e(g)T=T\hat{\rho}_v(g)$ for all $g$. Thus
$(\rho_e(g)T)^*=(T\hat{\rho}_v(g))^*$ and so
$T^*\rho_e(g)^*=(\hat{\rho}_v(g))^*T^*$, that is
$T^*\rho_e(g^{-1})=(\hat{\rho}_v(g^{-1}))T^*$. Restating this,
$T^*\rho_e(g)=(\hat{\rho}_v(g))T^*$, for all $g$. Thus,
suppressing some notation, $T^*T\hat{\rho}_v = T^*\rho_eT=
T^*T\hat{\rho}_v$. Since $T^*T$ commutes with $\hat{\rho}_v$ so
too does its square root $|T|$. We have $\rho_eU|T| =
U|T|\hat{\rho}_v = U\hat{\rho}_v|T|$ and it  follows, since $|T|$
has dense range for example, that
 $\rho_eU = U\hat{\rho}_v$ as desired.
\end{proof}

\subsubsection{Periodic frameworks.} We now show how the arguments of Section 3 can be applied to obtain  Fowler-Guest type formulae for  periodic bar-joint frameworks in $\bR^d$.
%\begin{equation}
% \Gamma_{p}(m) -\Gamma_{p}(s)=\Gamma_{p}(v)\cdot \Gamma(sp) - \Gamma_{p}(e)
% -\Gamma_{p}(rig).
%\end{equation}
The trace lists indicated in Theorem \ref{t:periodicFG}
are associated with finite-dimensional
representations of a finite group quotient $\G/\T$ of the spatial symmetry group $\G$, as we describe below.

Let $(G,p)$ be a countably infinite framework
in $\bR^d$ with distinct framework points and with spatial symmetry group $\G$ which contains a subgroup $\T$ isomorphic to $\bZ^d$ with $d$  independent generators $W_1, \dots, W_d$. It is in this sense that the framework is periodic.
We assume that the framework points are discrete in the sense that there are finitely many $\T-$orbits of framework points. With this condition
it follows  that $\G$ is a crystallographic group.
We do not assume that $\T$ is the minimal such subgroup.
In that case the quotient  $\G/ \T$
would be the associated point group of $\G$ but it is also of interest to consider
periodicity with respect to periods greater that the minimal period.

Consider the finite-dimensional Euclidean spaces $\H^{p}_v
\subseteq \H_v$ and
$\H^{p}_e \subseteq \H_e$ consisting of the vectors that are periodic with respect to $ \T$.  From the symmetry equations
\[
\rho_e(W_i)R(G,p) = R(G,p)\hat{\rho}_v(W_i), \quad i = 1, \dots , d,
\]
it follows readily that the rigidity matrix $R(G,p)$
determines a linear transformation
$R^{(p)}$ from  $H^{p}_v$ to $ \H^{p}_e$. The space $\ker R^{(p)}$
is the space $\ker R(G,p) \cap H^{p}_v$, which can be viewed as the space of periodic "infinitesimal" flexes for the framework $(G,p)$. Similarly the space
$\coker R^{(p)}$ is the space of periodic "infinitesimal" stresses.
(Note that a rotation flex $u=(u_v)_{v\in V}$, which is in the kernel of $R(G,p)$, is not a bounded sequence.)

The representations $\hat{\rho}_v, \rho_e$ of $\G$ induce representation
$\hat{\pi}_v, {\pi}_e$ of $\G /\T$ on the periodic vector spaces. Explicitly,
if $w=(w_f)_{f\in E}$ is in $\H^{p}_e$ then $\pi_e(g+\T)$ is well-defined by the
equation
\[
(\pi_e(g+\T)w)_f=w_{\sigma^{-1}(f)}.
\]
The representation $\hat{\pi}_v$ is defined similarly and the tensor
factorisation of $\hat{\pi}$  gives the  tensor factorisation
$\hat{\pi}_v=\pi_n \otimes \rho_{sp}$.

Since the rigidity matrix transformation  $R(G,p)$ and the transformations  $\hat{\rho}_v(h), \rho_e(h), h \in \G$,
leave invariant the spaces of periodic vectors we obtain from the
symmetry equations for $R(G,p)$ and $\hat{\rho}_v, \rho_e$ the induced symmetry equations
\begin{equation}
\pi_e(h)R^{(p)}  = R^{(p)}\hat{\pi}_v(h), \quad h \in \G/\T.
\end{equation}
As before, the representations $\hat{\pi},\pi_e $ do not depend on metrical detail and character lists for them are readily computable in terms of fixed elements.

Following the argument in Section 3, consider the orthogonal decompositions
\[
\H^{p}_v = \H^{p}_{v'} \oplus\H^{p}_{m}   \oplus \H^{p}_{rig}
\]
where $\H^{p}_{rig}= \H^{p}\cap \H_{rig}$ and $\H^{p}_{mech}$ is the complementary
space of $\H_{rig}^{per}$ in $\ker R^{(p)}$, and $\H^{p}_{v'}$ is the complementary
space of   $\ker R^{(p)}$  in $\H^{p}_v$. The rotational rigid motion flexes
are not periodic and so this intersection is a $d$-dimensional space corresponding to the translation flexes.
Similarly we have the decomposition
\[
\H^{p}_e=\H^{p}_{e'}\oplus \H^{p}_{str}.
\]

From the symmetry equations we see that  the component spaces
\[
\H^{p}_{v'},~~ \H^{p}_{m}, ~~\H^{p}_{rig}
\]
are invariant for $\hat{\pi}_v$ and so define subrepresentations of $\hat{\pi}_v$
whose trace lists we shall denote as
\[
\Gamma_{p}(v'), ~~\Gamma_{p}(m), ~~\Gamma_p({rig}),
\]
Similarly for the two subrepresentations of $\pi_e$ we obtain
the character lists
\[
 \Gamma_{p}(e'), ~~\Gamma_{p}(s)
\]
All five lists correspond to
some fixed suppressed set $h_1, \dots , h_s$  of generating elements of $\G/\T$.

\begin{thm}\label{t:periodicFG}
Let $(G,p)$ be a discrete periodic framework in $\bR^d$ with spatial symmetry group $\G$ and let  $\T \subseteq \G$ be a
full rank translation subgroup isomorphic to $\bZ^d$. Then
\begin{equation}
 \Gamma_{p}(m) -\Gamma_{p}(s)=\Gamma_{p}(v)\cdot \Gamma(sp) - \Gamma_{p}(e)
 -\Gamma_{p}(rig)
\end{equation}
where $ \Gamma_{p}(m)$ (resp. $\Gamma_{p}(s)$) are character lists for the representation
of the finite group  $\G/\T$ in the space of periodic (proper infinitesimal) mechanisms
(resp. the space of periodic stresses).
\end{thm}

\begin{proof}
The transformation $R^{(p)}$ induces an equivalence of the representations
which shows that $\Gamma_{p}(v')=\Gamma_{p}(e')$.
Since
\[
\Gamma(\hat{\pi_v})=\Gamma_p(v)\cdot \Gamma(sp)=\Gamma_{p}(v')+\Gamma_{p}(m)+\Gamma_p({rig})
\]
and
\[
\Gamma(\pi_e) = \Gamma(e)= \Gamma_{p}(e')+\Gamma_{p}(s)
\]
equation (11) follows.
\end{proof}

\begin{rem}
{\rm In the case of planar periodic frameworks evaluating at the
identity matrix give a periodic Maxwell rule, namely
\[
m_p-s_p = 2|V_p|-|E_p|-2
\]
where $m_p$ and $m_s$ are the dimension of the spaces of periodic
infinitesimal mechanisms and stresses, respectively, and $|V_p|$ and $|E_p|$
are the number of $\T$-orbits of vertices and edges, respectively.
In the periodic isostatic case $m_p=s_p=0$ (by definition) and we have the necessary
condition $2|V_p|-|E_p|-2$.

Periodic rigidity and isostaticity has been developed in interesting work of Ross\cite{ros} who has obtained a periodic version of Laman's theorem in the case that the
vertices in a unit cell for $\T$ are generically located.
We also note that Borcea and Streinu\cite{bor-str} have considered
more general forms of deformability of periodic frameworks. See also Owen and Power\cite{owe-pow-crystal}.
}
\end{rem}

%Another symmetry aspect which is novel in the infinite case is
%that a framework may be invariant under the dilation operation
%$(x,y) \to c(x,y)$ for some $c >0$.

\section{Symmetry in subframeworks and partitions}

We now show how latent symmetries can play a role in predicting
the singularity of asymmetric frameworks.

\subsection{Subframework symmetry}  Let $(G,p)$ be a proper bar-joint framework in
$\bR^2$ with a subframework $(X,p)$, where $X$ is a  subgraph of
$G$ (with at least one edge). Here, and below, it is convenient to use the redundant
notation $(X,p)$ with $p$ the full framework vector. The
Fowler-Guest formula holds for $(X,p)$ and in our notation takes
the form
\[
[\rho_{{mech}}^X]- [\rho_{st}^X]=  [\rho_{sp}^X]\cdot
[{\rho}_{n}^X] - [\rho_{e}^X]- [\rho_{rig}^X]
\]
where each $\rho^X$ is a representation of the spatial symmetry group of
$(X,p)$. In particular evaluating traces of the representations of
the identity symmetry gives the Calladine-Maxwell identity for
$(X,p)$, while evaluating at a reflection symmetry, $g$ say, gives
an identity which we write as
\[
m_X^g-s_X^g= 0 - b_X^g+1.
\]
Here  $b_X^g=\trace(\rho_{e}^X(g))$ is the number of framework edges  (bars) left
invariant by $g$. The
term $0$ arises from $\trace(\rho_{sp}^X(g))=0$, and for the
three-dimensional representation $\rho_{rig}^X$ we have
$\trace(\rho_{rig}^X(g))=-1$.

We now exploit the evident fact that the natural inclusion $\H^X_e
\subseteq \H^G_e$ respects stresses, that is, $\H^X_{st} \subseteq
\H^G_{st}$. This is simply because a vector in the cokernel of
$R(X,p)$ extends trivially to a vector in the cokernel of
$R(G,p)$.
%Suppose  that  $(G,p)$ is isostatic and that we can
%identify a subframework $(X,p)$ for which there is a reflection
%$g$ for which $- b_X^g+1$ is nonzero. Since $\H^G_{st}=0$ the
%space of stresses for $(X,p)$ is trivial and so $s_X^g=0$ and it
%follows that $m^g_X$ is nonzero, and $(X,p)$ has a nontrivial
%space of flexes. However, it may be evident, as in our example,
%that any flex of $(X,p)$ has an extension to a flex of $(G,p)$ in
%which case we have a contradiction and we infer that $(G,p)$ is
%not isostatic.
The following theorem  gives
a family of necessary conditions all of which are computable by
simple counting.

Combining these facts we obtain

\begin{thm}\label{t:partitiontheorem}
Let $(G,p)$ be a proper isostatic framework in $\bR^d$. Then

(i) for
each proper subframework $(X,p)$ and each spatial symmetry $g$ of
$(X,p)$ we have
\[
|\trace(g).v_X^g - e_X^g - \trace(\rho_{rig}(g))| \leq dv_X-e_X
-d(d+1)/2
\]
{where $v_X^g$ (resp. $e_X^g$) is the number of vertices (resp.
edges) in the graph $X$ that are unmoved by the symmetry.}

(ii) For planar frameworks  a necessary condition for
isostaticity is that for each reflection symmetry $g$ of a
subframework $(X,p)$
\[
|-e_X^g+1|\leq 2v_X-e_X-3.
\]
\end{thm}

\begin{proof}
In  $(X,p)$ we have
\[
m_X=dv_X-e_X-d(d+1)/2,
\]
which follows on evaluating the general formula at the identity
symmetry and noting as above that $s_X=0$. For the symmetry $g$ of $(X,p)$
 we have $|m_X^g| \leq m_X$,
since $m_X$ is the dimension of the mechanism space of $(X,p)$. On
the other hand the evaluation of  traces on the identity element  gives
\[
m_X^g - s_X^g =\trace(g).v_X^g - e_X^g - \trace(\rho_{rig}(g))
\]
Combining these facts we obtain (i), from which (ii) follows.
\end{proof}

\vspace{.6in}
\begin{figure}[h]
\centering
\includegraphics[width=5cm]{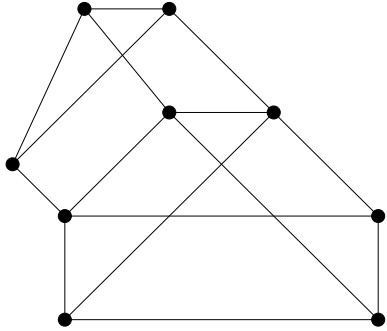}
\caption{A framework with reflection symmetry in a sub-graph and a
singular Jacobian. }
\end{figure}

The second part of the theorem is illustrated in Figure 4 where there is an evident
subframework $X$ with six vertices with a mirror symmetry. Since the inequality of the theorem is violated for $X$ the entire framework fails to be isostatic.

\subsection{Partition symmetry} We now
show how symmetries associated with vertex partitioning can be
significant for singularity. The idea here is that on removing the
framework edges connecting vertices within each of the sets of a
partition of $V$ one may be left with a set of "crossing" edges
which has evident symmetry. In this event one can add edges to
create complete graph frameworks within the partition sets thereby
creating  a body framework. If, by  symmetry and counting
conditions, the resulting framework  has proper flexes then the original
framework inherits the same proper flexes. This situation occurs
for example in the simple framework of Figure 5.

More precisely let $G=(V,p)$ be a framework in $\bR^2$, where each
vertex has degree greater than 1, and let $V_1, \dots ,V_n$ be a
partition of $V$. Let
\[
\S =\{V_1, \dots ,V_n, e_1,\dots ,e_m\}
\]
where $e_1,\dots ,e_m$ are the edges of $G$ which have vertices in
distinct partition sets. Delete from $p$ the framework points
which are not endpoints of the edges $e_i$ to create a framework
vector $p'$ (representing pins). Then $(\S, p')$ is a body
framework and we say that it is derived from $(G,p)$, or that it
is a partition-derived body framework. Note that for  a trivially
derived body framework, where each partition set is a singleton,
the total number of point body occurrences is the sum of the
degrees of the vertices in $G$, which is $2e$. Thus $c = 2e$ and
the isostatic condition in the trivially derived framework gives
$2c = 2n+3(e-1)$, which implies $e = 2n-3$ as  expected.

The following theorem, together with Theorem \ref{t:body} give
necessary conditions for isostaticity.

\begin{thm}\label{t:last}
Let $(G,p)$ be a framework in $\bR^2$ and let $(\S,p)$ be a
partition-derived body bar framework. Then

(i) a (non-trivial) flex of $(\S, p)$ gives a (non-trivial) flex of
$(G,p)$.

(ii) if $(G,p)$ is isostatic then a reflection symmetry of $(\S,p)$
fixes exactly one edge of $(\S,p)$.
\end{thm}

\begin{proof} Let the set of velocity vectors
$\{u_i, v_e, a_e\}$ be a flex of $(\S,p)$. For any two points
$p_i$ and $p_j$ in body $e$,  $u_i = v_e + a_e(p_i-p_e)^{\pi/2},
u_j = v_e + a_e(p_j-p_e)$. Thus $u_i - u_j = a_e(pi-pj)^{\pi/2}$
and $(u_i - u_j ).(p_i-p_j) = 0$. Since every pair of points
joined by a framework edge are both in some body of $\S$ it
follows that the set $\{u_i\}$ is a flex of $(G,p)$. Now (i)
follows and (ii) follows from (i).
\end{proof}

\vspace{.5in}

\begin{center}
\begin{figure}[h]
\centering
\includegraphics[width=4cm]{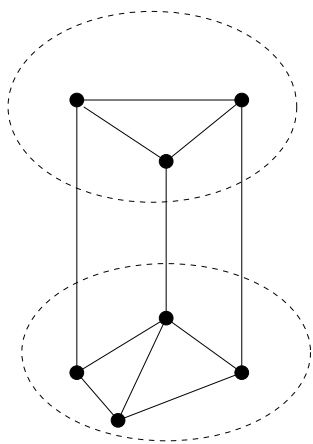}
\caption{A framework with vertical reflection symmetry in a
partition derived graph and a singular Jacobian. }
%The ellipses
%indicate the vertex partition.}
\end{figure}
\end{center}

\end{document}